\documentclass[reqno]{amsart}
\usepackage{amssymb,eucal}

\newcommand{\setmin}{{\setminus}}

\newcommand{\jsd}{join-sem\-i\-dis\-trib\-u\-tive}

\newcommand{\set}[1]{\{\,#1\,\}}
\newcommand{\setm}[2]{\{\,#1\mid#2\,\}}

\newcommand{\BB}[1]{\mathbf{\mathcal{B}_{#1}}}

\newcommand{\ovY}{\overline{Y}}
\newcommand{\R}{\mathbb{R}}

\DeclareMathOperator{\At}{At}

\DeclareMathOperator{\Co}{Co}
\DeclareMathOperator{\Ex}{Ex}

\newcommand{\interior}[1]{int_{#1}}
\newcommand{\DD}{\mathbin{D}}

\newcommand{\Subm}[1]{\mathrm{Sub}_{\wedge}\mathbf{\mathcal{B}_{#1}}}
\newcommand{\CoX}{\Co(\mathbb{R}^n,X)}

\theoremstyle{plain}

\newtheorem{lemma}{Lemma}[section]
\newtheorem{theorem}[lemma]{Theorem}
\newtheorem{proposition}[lemma]{Proposition}

\newtheorem{claim}{Claim}

\newtheorem*{stat}{\name}
\newcommand{\name}{testing}

\theoremstyle{definition}

\newtheorem{example}[lemma]{Example}
\newtheorem{problem}{Problem}

\theoremstyle{remark}

\newcommand{\qedc}{{\qed}~{\rm Claim~{\theclaim}.}}

\newenvironment{cproof}
{\begin{proof}[Proof of Claim.]}
{\qedc\renewcommand{\qed}{}\end{proof}}

\begin{document}

\author{K. V.~Adaricheva}
\address{Institute of Mathematics of SB RAS, Acad. Koptyug Prosp., 4,
630090, Novosibirsk, Russia}

\email{ki13ra@yahoo.com}
\keywords{Lattice, relatively convex set, join-semidistributive, lower bounded,
convex geometry. }
\subjclass{Primary 06B15, 51E99, 52A20. Secondary: 51D20, 05B25}

\title{Join-semidistributive lattices of relatively convex sets}

\begin{abstract}
We give two sufficient conditions for the
lattice $\CoX$ of
relatively convex sets of $\R^n$ to be \jsd, where $X$ is a finite
union of segments. We also prove that every finite lower bounded
lattice can be embedded into
$\CoX$, for a suitable finite subset
$X$ of $\R^n$.
\end{abstract}

\maketitle

\section{Introduction}
A lattice $L$ is \emph{\jsd}, if
\[
 x \vee y = x \vee z \text{ implies that }
 x \vee y = x \vee (y\wedge z),
\]
for all $x,y,z\in L$.
Let $X\subseteq\R^n$, and let
$\CoX$ denote the lattice of convex subsets of $\R^n$
\emph{relative to} $X$, that is,
\[
\CoX = \setm{Y\subseteq\R^n}{Y=\Co(Y)\cap X},
\]
where $\Co(Y)$ denotes \emph{the convex hull} of $Y$,
for any $Y \subseteq \R^n$.
For all $X\subseteq\R^n$, the closure operator $\phi\colon\BB{X}\to\BB{X}$,
where $\phi(Y)=\Co(Y)\cap X$ for all
$Y\subseteq\R^n$, satisfies the so-called \emph{anti-exchange axiom}
that makes lattices of relatively convex sets just another example
of \emph{a convex geometry} (see
the extensive monograph~\cite{KLS}, also~\cite{AGT}).
It is well known (cf.~\cite{AGT})
that a finite convex geometry is \jsd, whence the lattice $\CoX$ is \jsd, for
any finite $X\subseteq\R^n$.

Problem 3 in~\cite{AGT} asks about a description of lattices embeddable
into lattices of the form $\CoX$ with finite $X$.
Since any sublattice of a \jsd\ lattice is \jsd\ itself, all those lattices
must also be \jsd.
Although the current paper does not provide a solution of the problem, it
suggests some approaches to it. The main idea is to consider a more
general setting for the problem dropping the
requirement for $X$ to be finite.

For a lattice $L$ with the least element $0_L$, let
$\At(L)$ denote the set of \emph{atoms} of $L$, that is,
$\At(L)=\setm{x\in L}{0_L\prec x}$.
While finite convex geometries are always \jsd,
a convex geometry $L$ satisfies a weaker property:
\[
 x \vee y = x \vee z \text{ implies that }
 x \vee y = x \vee (y\wedge z),
\]
for all $x\in L$ and all $y,z\in\At(L)$. In other words, if
$x\vee y=x\vee z$, for some $x\in L$ and $y,z\in\At(L)$ the either
$y=z$ or $y,z\leq x$.
How weak this property is can be seen from the following result established
in~\cite{B}: \emph{every finite lattice can be embedded into
$\CoX$, for some $n\in\omega$ and $X\subseteq\R^n$}.
Thus we would like to generalize Problem~3 from~\cite{AGT}, dropping
the requirement for $X$ to be finite but still assuming $\CoX$ to be \jsd:

\begin{problem}
Which finite lattices can be embedded into \jsd\ lattices of the
form $\CoX$?
\end{problem}

It turns out that sets $X$ for which the corresponding lattice $\CoX$ is
\jsd\ are quite specific. The third section of the paper is mostly devoted
to the case when $X$ is a finite union of segments, which seems to be a
natural generalization of finiteness of $X$. We provide two sufficient
conditions for $X$ to ensure $\CoX$ to be \jsd.

The last section is devoted to an important proper subclass of the class
of \jsd\ lattices, the class of so-called \emph{lower bounded
lattices}. We prove that every finite lower bounded lattice embeds
into a finite lower bounded lattice of the form $\CoX$.
Another proof of this result can be found also in~\cite{WS}.

Here we use an essentially geometric idea, first constructing an embedding
of the lattice $\Subm{n+1}$
of meet-subsemilattices of the Boolean lattice $\BB{n+1}$ into
the lattice of bounded convex subsets of $\R^n$, and then finding
a finite set $X$ which provides an embedding into $\CoX$.
We hope that this construction might give some additional insight into the
question whether every finite \jsd\ lattice embeds
into a finite lattice $\CoX$.

\section{Basic concepts}
For any $a,b\in \R^n$, let $(a,b)$ denote the open segment
and let $[a,b]$ denote the closed segment whose end points are $a$ and $b$,
that is,
\begin{align*}
(a,b)&=\setm{x\in\R^n}{x=\lambda a+(1-\lambda)b\text{ for some }\lambda\in(0,1)},\\
[a,b]&=\setm{x\in\R^n}{x=\lambda a+(1-\lambda)b\text{ for some }\lambda\in[0,1]}.
\end{align*}
It is straightforward to verify that for any
$Y\subseteq\R^n$,
\[
\Co(Y)=\bigcup_{i\in\omega}Y^{(i)},
\]
where $Y^{(0)}=Y$ and $Y^{(i+1)}=\setm{[a,b]}{a,b\in Y^{(i)}}$, for all
$i\in\omega$.

A convex subset $F\subseteq P$ of a convex potytope $P$ is \emph{a face of} $P$, if
$(a,b)\cap F\ne\varnothing$ implies $[a,b]\subseteq F$,
for all $a,b\in P$. An element $x$ of a convex set $X\subseteq\R^n$ is
\emph{an extreme point} of $X$ if $x\notin\Co(X\setmin\set{x})$. Let
$\Ex(X)$ denote the set of extreme points of $X$, for any $X\in\Co(\R^n)$.

For any $Y\subseteq\R^n$, we denote by $\ovY$ the closure
of $Y$ and by $\interior{n}(Y)$ the interior of $Y$ in the Euclidean
topology of $\R^n$.

\begin{lemma}\label{vertex}
Let $X\subseteq \R^n$ be a finite union of segments.
Then $\Co(\overline{X})=\overline{\Co(X)}$.
In particular, if $x\in\Ex(\overline{\Co(X)})$
then $x$ is an extreme point of a closure of a segment from $X$.
\end{lemma}

\begin{proof}
The proof is straightforward.
\end{proof}

\begin{lemma}\label{hyper}
Let $P\subseteq\R^n$ be a convex polytope and let $F$ be a face of $P$.
Then $\Co(Y)\cap F=\Co(Y\cap F)$, for any $Y\subseteq P$.
\end{lemma}

\begin{proof}
By induction on $k$, we prove that $Y^{(k)}\cap F\subseteq(Y\cap F)^{(k)}$,
for all $k\in\omega$. For $k=0$, the conclusion is obvious. Let $k>0$ and let
$x\in Y^{(k)}\cap F$. Then there exist $a,b\in Y^{(k-1)}$ such that
$x\in[a,b]$. If $x=a$ or $x=b$, then
$x\in Y^{(k-1)}\cap F\subseteq(Y\cap F)^{(k-1)}$ by the
induction hypothesis. Otherwise, $x\in(a,b)\cap F$,
whence $a,b\in F$ since $F$ is a face of $P$. Therefore,
$a,b\in Y^{(k-1)}\cap F\subseteq(Y\cap F)^{(k-1)}$ by the induction hypothesis,
whence $x\in(Y\cap F)^{(k)}$.
\end{proof}

For any $Y\subseteq\R^n$, let $\psi_Y\colon\Co(\R^n)\to\Co(\R^n,Y)$ be the
map defined by $\psi_Y(X)=X\cap Y$, for any $X\in\Co(\R^n)$. Then $\psi_Y$
preserves meets, for any $Y\subseteq\R^n$.

\begin{lemma}\label{hom}
Let $P$ be a convex polytope and let $X\subseteq P$. Then
the map $\psi_F\colon\CoX\to\Co(\R^n,X\cap F)$ defined by
$\psi_F(Y)=Y\cap F$ is a surjective lattice homomorphism, for any
face $F$ of $P$.
\end{lemma}

\begin{proof}
The surjectivity of $\psi_F$ follows from the fact that if
$A=\Co(A)\cap X\cap F$ then $A=\psi_F(\Co(A)\cap X)$.
Let $A,B\in\CoX$. Evidently, $\psi_F$ preserves meets.
Applying Lemma~\ref{hyper} we get
\begin{align*}
\psi_F(A\vee B)&=\Co(A\cup B)\cap X\cap F=
\Co\bigl((A\cap F)\cup(B\cap F)\bigr)\cap X\\
&=\bigl(\Co(A\cap F)\cap X\bigr)\vee\bigl(\Co(B\cap F)\cap X\bigr)=
\psi_F(A)\vee\psi_F(B),\\
\end{align*}
whence $\psi_F$ preserves joins.
\end{proof}

\section{Join-semidistributivity of $\CoX$}
If $X \subseteq \R^n$ is finite, then, as we mentioned above, the lattice $\CoX$ is a finite convex
geometry; in particular, it is \jsd. However, we do not know how far this
fact can be extended.

\begin{problem}\label{describeX}
Describe sets $X \subseteq \R^n$ such that the lattice $\CoX$
is \jsd.
\end{problem}

To remind that not every $X$ suits, we recall an example
given in~\cite{B}.

\begin{example}\label{triangle}
Let $X$ contain the ($2$-dimensional) interior of some triangle $TML$.
Pick any point
$K$ inside that interior.
Then the interior of each triangle $TMK$, $TLK$, and $MLK$ belongs to
$\CoX$, and they form a modular sublattice isomorphic to $M_3$.
In particular, $\CoX$ is not \jsd.
\end{example}

A subset $X$ of $\R^n$ is \emph{sparse}, if
$\interior{2}(X\cap H)=\varnothing$, for any
$2$-dimensional affine subspace $H$ of $\R^n$.
From Example~\ref{triangle}, it follows that every
set $X$ satisfying the requirement of Problem~\ref{describeX} has to be
sparse.

Observe that if $X$ is a line in $\R^n$ then $\CoX$ is isomorphic to $\Co(\R)$,
the lattice of order convex subsets of $\R$, and the latter is \jsd\
(see Theorem~14 in~\cite{BB}).

Another extreme case is when $X$ is the boundary
of a ball; in this case, the lattice $\CoX$
is Boolean (cf.~an example of section~9 in~\cite{B}); in particular, it is
distributive.
This gives two natural examples of sparse sets which qualify for
Problem~\ref{describeX}.
Unfortunately, being a sparse set is a necessary condition
but not sufficient.

\begin{example}\label{3-lines}
Let $X$ be the union of three lines $A$, $B$, and $C$ which are on the same
plane and have a common intersection. Then
$A \vee B = A \vee C = X$ but $A \vee (B \cap C) = A$ in $\CoX$.
\end{example}

On the other hand, if we take \emph{segments} instead of lines,
then the corresponding lattice turns out to be \jsd. Thus the following
question is rather natural: \emph{if} $X$ \emph{is a finite union of segments,
is the lattice} $\CoX$ \emph{\jsd?}
Unfortunately, even this simplest generalization of finiteness of $X$
does not ensure that $\CoX$ is \jsd, as the example below demonstrates.

\begin{example}\label{PM}
Let $T$ be a triangle in $\R^2$ with the set of extreme points $\set{a,b,c}$
and let $p,m \in \interior{2}T$, $p\ne m$. Without loss of generality, we
may assume that
$p$, $m$, and $a$ are not collinear. We put
$X=[b,c]\cup[p,a]\cup[m,a]$ and $A = [b,c]$, $B = (p,a)$, $C =(m,a)$.
Then
$A\vee B=A\vee C=X\setmin\set{a}\ne A\vee(B\wedge C)=A$ in
$\Co(\R^2,X)$.
Thus this lattice is not \jsd.
\end{example}

We note that the failure of join-semidistributivity in the example above is due to
the fact that closed segments $[p,a]$ and $[m,a]$ have a common point.
Also, it is essential that $(p,a)$ and
$(m,a)$ are subsetes of $\interior{2}T$.
Were points $p$ and $m$ chosen, say, on
faces $[a,b]$ and $[a,c]$ of the triangle $T$, respectively,
the lattice $\CoX$ would be \jsd.

For the rest of this section, we assume $X$ to be a finite union of segments.
The following theorem provides two sufficient conditions for $\CoX$ to be \jsd.
Each of them eliminates at least one condition
that plays role in Example~\ref{PM}.

\begin{theorem}
Let $n,k\in\omega$ and let $X=\bigcup\setm{I_j}{j<k}$, where
$I_j\subseteq\R^n$ is a segment, for all $j<k$.
Consider the following two conditions:
\begin{enumerate}
\item
$\overline{I_s}\cap\overline{I_t}=\varnothing$, for all $s,t<k$, $s\ne t$;
\item
there exists a convex polytope $P\subseteq\R^n$ such that
for any $j<k$, $I_j$ is a subset of a face of $P$.
\end{enumerate}
If $X$ satisfies either $\mathrm{(i)}$ or $\mathrm{(ii)}$ then the lattice
$\CoX$ is \jsd.
\end{theorem}

\begin{proof}
We agrue by induction on $n$. Let $n=1$. For any $X\subseteq\R$, the lattice
$\Co(\R,X)$ is the lattice of order-convex subsets of $X$ endowed with the
standard (linear) order, thus it is \jsd\ (see~\cite[Theorem~14]{BB}).

Let $n>1$.
Suppose that $X$ satisfies either (i) or (ii) and
$A\vee B=A\vee C>A\vee(B\cap C)$, for some $A,B,C\in\CoX$.
Let $Y=\Co(A\vee(B\cap C))$.
Then $B,C\not\subseteq Y$. We prove that there are a convex polytope $Q$
and a face $F$ of $Q$ such that $B\cap F\not\subseteq Y$ and
$Y\subseteq Q$.

Suppose first that $X$ satisfies (i). By Lemma~\ref{vertex}, we get
\[
K=\overline{\Co(A\cup B)}=\Co(\overline{A\vee B})=\Co(\overline{A\vee C})=
\overline{\Co(A\cup C)}.
\]
If $K\not\subseteq\ovY$, then there exists an extreme point $a\in\Ex(K)$
such that
$a\notin\ovY$. Since $\overline{A}\subseteq\ovY$, by Lemma~\ref{vertex},
$a\in\overline{B}\cap\overline{C}$ contradicting (i). Thus,
$B\subseteq K\subseteq\ovY$ but
$B\not\subseteq Y$. Therefore, there exists a face $F$ of $\ovY$ such
that $B\cap F\not\subseteq Y$. We take $Q=\ovY$ in this case.

Suppose that $X$ satisfies (ii). Since $B\not\subseteq Y$, there is a face
$F$ of $P$ such that $B\cap F\not\subseteq Y$. We take $Q=P$ in this case.

By Lemma~\ref{hom}, the map
$\psi_F\colon\Co(\R^n,X\cap Q)\to\Co(\R^n,X\cap Q \cap F)$ is a lattice homomorphism.
Thus,
$\psi_F(A)\vee\psi_F(B)=\psi_F(A)\vee\psi_F(C)$. Also, the lattice
$\Co(\R^n,X\cap F)$ is isomorphic to the lattice
$\Co(\R^m,X\cap F)$,
where $m\in\omega$ is the dimension of an affine subspace of $\R^n$
containing $F$. Moreover, $X\cap F$ is a finite union of segments.
By the induction hypothesis, the lattice $\Co(\R^m,X\cap F)$ is \jsd, whence
\begin{align*}
B\cap F=&\psi_F(B)\subseteq\psi_F(A\vee B)=\\
&\psi_F(A)\vee\bigl((\psi_F(B)\cap\psi_F(C))\bigr)=\\
&\psi_F\bigl(A\vee(B\cap C)\bigr)=\psi_F(Y)\subseteq Y,
\end{align*}
a contradiction.
\end{proof}

\section{Lower bounded lattices as sublattices of finite $\CoX$}
In this section, we consider sublattices of lattices of the form $\CoX$,
where $X\subseteq \R^n$ is finite.
As was observed in~\cite{AGT}, we do not know yet any special type of
finite convex geometries which admit any finite \jsd\ lattice as a sublattice. 
We have a partial confirmation that lattices of the form $\CoX$ could be
such a  "universal'' class of convex geometries for the class of finite
\jsd\ lattices.

The main result of this section shows that, at least, this class is universal
for the class of
finite \emph{lower bounded lattices} which is a proper subclass in the class
of finite \jsd\ lattices.
We recall that a (finite) lattice is \emph{lower bounded}, if it is an
image of a finitely generated free lattice under
\emph{a lower bounded homomorphism}, that is,
the preimage of every element under this homomorphism has a least element.
We refer the reader to the comprehensive monograph on the topic~\cite{FJN}.
There exist at least two other particular classes of finite convex geometries
which admit every finite lower
bounded lattice as a sublattice: suborder lattices of finite
partial orders~\cite{Si} and
subsemilattice lattices of finite semilattices~\cite{Ad,Rep}.

Unlike these known examples, lattices of relatively convex subsets are
\emph{not} necessarily lower bounded. The simplest example is
$\Co(\R,X)$, where $X$
consists of four different points on the same line.
The other common feature of many types of convex geometries is that
they are biatomic. Due to~\cite{BB}, a lattice $L$ with the least element
$0_L$ is \emph{biatomic} if
for any $x\in\At(L)$ and any $y,z\in\At(L)$,
the inequality $x\leq y\vee z$ implies that there are
$y',z'\in\At(L)$ such that $y'\leq y$, $z'\leq z$, and $x\leq y'\vee z'$.

A result from~\cite{AW} shows that \emph{not} every finite \jsd\ lattice
embeds
into a finite biatomic \jsd\ lattice. The counter-example from~\cite{AW}
is the lattice $\Co(\R^2,X)$, where $X$ is a $5$-element set of points
on a plane. In particular, this emphasizes that lattices of relatively
convex subsets are
essentially non-biatomic, thus might serve as a ``universal'' class of convex
geometries for the class of finite \jsd\ lattices.

Observe that an alternate approach which leads to the result that
every finite lower bounded lattice is a sublattice of some $\CoX$ with
finite $X$ is presented in~\cite{WS}.
The authors of~\cite{WS} find an embedding of every finite lower bounded
lattice into the lattice of convex polytopes of a finite-dimensional
vector space, from where the result easily follows.

\begin{proposition}\label{step1}
For every $n<\omega$, the lattice $\Subm{n+1}$ embeds into the lattice
of bounded convex sets of $\R^n$.
\end{proposition}

\begin{proof}
Let $S_{\mathbf{n+1}}$ denote a regular polytope in $\R^n$ with $n+1$ vertices.
It is not that important to have a \emph{regular} polytope,
but it is easier to deal with because of the total symmetry of the argument.
Thus, in $\R^2$ it is an equilateral triangle,
in $\R^3$ it is a regular tetrahedron, etc.

Let $\Ex(S_{\mathbf{n+1}})=\setm{p_i}{i\leqslant n+1}$.
We define the map $\psi\colon\BB{n+1}\to\Co(\R^n)$ by the rule
\begin{equation}\label{psi}
\psi(t)=\begin{cases}
\varnothing,\text{ if }t=\mathbf{n+1},\\
\set{p_i},\text{ if }\mathbf{n+1}\setmin t=\set{i},\\
\interior{|A|}\Co\bigl(\setm{p_i}{i\in A =\mathbf{n+1}\setmin t}\bigr),
\text{ if }|t|<n.
\end{cases}
\end{equation}

\begin{claim}\label{join}
For any $a,b\in\BB{n+1}$,
$\Co\bigl(\psi(a)\cup\psi(b)\bigr)=\psi(a)\cup\psi(b)\cup\psi(a\cap b)$.
\end{claim}

\begin{cproof}
Without loss of generality, we may assume that
$a$ and $b$ are noncomparable.
By induction on $i$, we prove that
$\bigl(\psi(a)\cup\psi(b)\bigr)^{(i)}
\subseteq\psi(a)\cup\psi(b)\cup\psi(a\cap b)$, for all $i\in\omega$.
For $i=0$, the conclusion is obvious. Suppose that
$i<\omega$ and that
$z\in\bigl(\psi(a)\cup\psi(b)\bigr)^{(i+1)}
\setmin\bigl(\psi(a)\cup\psi(b)\bigr)^{(i)}$.
Then there are $\lambda\in(0,1)$, $x,y\in\bigl(\psi(a)\cup\psi(b)\bigr)^{(i)}$
such that $z=\lambda x+(1-\lambda)y$. By the induction hypothesis,
$x,y\in\psi(a)\cup\psi(b)\cup\psi(a\cap b)$.
We consider several cases:

\textbf{Case~1.}
$x,y\in\psi(a)$ or $x,y\in\psi(b)$. In this case,
$z\in\psi(a)\cup\psi(b)$ since both $\psi(a)$ and $\psi(b)$ are convex.

\textbf{Case~2.}
$x\in\psi(a)$ and $y\in\psi(b)$. In this case, there are $\lambda_k\in(0,1)$,
$k\in\mathbf{n+1}\setmin a$, and $\mu_l\in(0,1)$,
$l\in\mathbf{n+1}\setmin b$, such that
\begin{align*}
&\sum\setm{\lambda_k}{k\in\mathbf{n+1}\setmin a}=
\sum\setm{\mu_l}{l\in\mathbf{n+1}\setmin b}=1\text{ and}\\
&x=\sum\setm{\lambda_k p_k}{k\in\mathbf{n+1}\setmin a},
\qquad y=\sum\setm{\mu_l p_l}{l\in\mathbf{n+1}\setmin b}.
\end{align*}
Then
\[
z=\sum\setm{\lambda\lambda_k p_k}{k\in\mathbf{n+1}\setmin a}+
\sum\setm{(1-\lambda)\mu_l p_l}{l\in\mathbf{n+1}\setmin b}.
\]
Moreover, $\lambda\lambda_k$, $(1-\lambda)\mu_l\in(0,1)$, for all
$k\in\mathbf{n+1}\setmin a$ and all $l\in\mathbf{n+1}\setmin b$, and
\[
\sum\setm{\lambda\lambda_k}{k\in\mathbf{n+1}\setmin a}+
\sum\setm{(1-\lambda)\mu_l}{l\in\mathbf{n+1}\setmin b}=
\lambda\cdot 1+(1-\lambda)\cdot 1=1.
\]
Thus, $z\in\psi(a\cap b)$.

\textbf{Case~3.}
$x\in\psi(a)$, $y\in\psi(a\cap b)$. In this case, there are
$\lambda_k\in(0,1)$,
$k\in\mathbf{n+1}\setmin a$, and $\mu_l\in(0,1)$,
$l\in\mathbf{n+1}\setmin(a\cap b)$, such that
\begin{align*}
&\sum\setm{\lambda_k}{k\in\mathbf{n+1}\setmin a}=
\sum\setm{\mu_l}{l\in\mathbf{n+1}\setmin(a\cap b)}=1\text{ and}\\
&x=\sum\setm{\lambda_k p_k}{k\in\mathbf{n+1}\setmin a},
\qquad y=\sum\setm{\mu_l p_l}{l\in\mathbf{n+1}\setmin(a\cap b)}.
\end{align*}
Then
\[
z=\sum
\setm{\bigl(\lambda\lambda_k+(1-\lambda)\mu_k\bigr)p_k}{k\in\mathbf{n+1}\setmin a}+
\sum\setm{(1-\lambda)\mu_l p_l}{l\in a\setmin b}.
\]
Again, all the coefficients are from $(0,1)$, and
\begin{align*}
\sum&\setm{\lambda\lambda_k+(1-\lambda)\mu_k}{k\in\mathbf{n+1}\setmin a}+
\sum\setm{(1-\lambda)\mu_l}{l\in a\setmin b}=\\
=&\lambda\sum\setm{\lambda_k}{k\in\mathbf{n+1}\setmin a}+
(1-\lambda)\sum\setm{\mu_l}{l\in\mathbf{n+1}\setmin(a\cap b)}=\\
=&\lambda\cdot 1+(1-\lambda)\cdot 1=1.
\end{align*}

Thus, $z\in\psi(a\cap b)$.
Therefore, we have proved that
$\Co\bigl(\psi(a)\cup\psi(b)\bigr)\subseteq
\psi(a)\cup\psi(b)\cup\psi(a\cap b)$.

We prove the inverse inclusion. It suffices to show that
$\psi(a\cap b)\subseteq\Co\bigl(\psi(a)\cup\psi(b)\bigr)$. Let
$z\in\psi(a\cap b)$. There are $\lambda_k\in(0,1)$,
$k\in\mathbf{n+1}\setmin(a\cap b)$
such that $\sum\setm{\lambda_k}{k\in\mathbf{n+1}\setmin(a\cap b)}=1$ and
\[
z=\sum\setm{\lambda_k p_k}{k\in\mathbf{n+1}\setmin(a\cap b)}.
\]
We put
\begin{align*}
\lambda&=\Bigl(\sum\setm{\lambda_k}{k\in b\setmin a}+
\frac{1}{2}\sum\setm{\lambda_k}{k\in\mathbf{n+1}\setmin(a\cup b)}\Bigr)^{-1};\\
x&=\sum\setm{\frac{\lambda_k}{\lambda}p_k}{k\in b\setmin a}+
\sum\setm{\frac{\lambda_k}{2\lambda}p_k}{k\in\mathbf{n+1}\setmin(a\cup b)};\\
y&=\sum\setm{\frac{\lambda_k}{1-\lambda}p_k}{k\in a\setmin b}+
\sum\setm{\frac{\lambda_k}{2(1-\lambda)}p_k}{k\in\mathbf{n+1}\setmin(a\cup b)}.
\end{align*}
We get
\begin{align*}
\sum&\setm{\frac{\lambda_k}{\lambda}}{k\in b\setmin a}+
\sum\setm{\frac{\lambda_k}{2\lambda}}{k\in\mathbf{n+1}\setmin(a\cup b)}=\\
=&\frac{1}{\lambda}\Bigl(\sum\setm{\lambda_k}{k\in b\setmin a}+
\frac{1}{2}\sum\setm{\lambda_k}{k\in\mathbf{n+1}\setmin(a\cup b)}\Bigr)=\\
=&\frac{1}{\lambda}\cdot\lambda=1;\\
\sum&\setm{\frac{\lambda_k}{1-\lambda}}{k\in a\setmin b}+
\sum\setm{\frac{\lambda_k}{2(1-\lambda)}}{k\in\mathbf{n+1}\setmin(a\cup b)}=\\
=&\frac{1}{1-\lambda}\Bigl(\sum\setm{\lambda_k}{k\in a\setmin b}+
\frac{1}{2}\sum\setm{\lambda_k}{k\in\mathbf{n+1}\setmin(a\cup b)}\Bigr)=\\
=&\frac{1}{1-\lambda}\cdot(1-\lambda)=1.\\
\end{align*}
Thus, $x\in\psi(a)$ and $y\in\psi(b)$.
Moreover, $z=\lambda x+(1-\lambda)y$,
whence $z\in\Co\bigl(\psi(a)\cup\psi(b)\bigr)$.
\end{cproof}

For any $S\in\Subm{n+1}$, we put
\begin{equation}\label{psi'}
\varphi(S)=\bigcup\setm{\psi(t)}{t\in S}.
\end{equation}

According to Claim~\ref{join}, $\varphi(S)\in\Co(\R^n)$, for any
$S\in\Subm{n+1}$.
We verify that $\varphi$ is a lattice homomorphism
from $\Subm{n+1}$ to $\Co(\R^n)$.
It is straighforward that $\varphi$ is one-to-one. Moreover, $\varphi$
preserves meets.

Let $S_0,S_1\in\Subm{n+1}$ and let
$S=S_1\vee S_2$.
If $t\in S\setmin(S_0\cup S_1)$,
then $t=t_0\cap t_1$, for some $t_i\in S_i$, $i<2$. Hence, by
Claim~\ref{join},
$\psi(t)\subseteq\Co\bigl(\psi(t_0)\cup\psi(t_1)\bigr)\subseteq\varphi(S_0)
\vee\varphi(S_1)$. Thus
$\varphi(S_0\vee S_1)\subseteq\varphi(S_0)\vee\varphi(S_1)$, whence
$\varphi$ preserves joins.
\end{proof}

For any $k<\omega$, for any $\lambda\geqslant 0$ small enough,
and for any convex
polytope $P\subseteq\R^k$, let $P^\lambda$ denote the (nonempty) convex polytope
which is a subset of $P$, whose faces are parallel to the corresponding
faces of $P$, and
$\rho(P^\lambda,P)=\lambda$, where $\rho(A,B)$ denotes the distance between
$A$ and $B$ defined by the standard Euclidean metric $\rho$.
For any $x\in\Ex P$, let $x^\lambda$ denote the
corresponding extreme point of $P^\lambda$.

We fix $n\in\omega$ and consider the polytope $S_{\mathbf{n+1}}$ defined in the
proof of Proposition~\ref{step1}. Let $\lambda>0$ be small enough.

If $A\subseteq\mathbf{n+1}$ and $|A|=k+1$, for some $k<\omega$, then
$S_A$ denotes the regular polytope in $\R^k$ with the set of extreme points
$\Ex S_A=\setm{p_i}{i\in A}$. For any
$B\subseteq A$, we put
\[
H_B=\setm{\sum_{i\in B}\lambda_i p_i^\lambda}{\lambda_i\in\R
\text{ for all }i\in B}.
\]
For any different $i,j\in A$, let $p(i,A,j)$ be a unique point from the
intersection $[p_i,p_j]\cap H_{A\setmin\set{j}}$. We put
\[
T(A,\lambda,j)=\Co\bigl(\setm{p_i,p(i,A,j)}{i\in A,i\ne j}\bigr).
\]
For any $j\in A$, the convex polytope $T(A,\lambda,j)$ has two parallel faces:
one is the face $S_{A\setmin\set{j}}$ of the polytope
$S_A$, the other is the face
$S'_{A\setmin\set{j}}=\Co\bigl(\setm{p(i,A,j)}{i\in A,i\ne j}\bigr)$.

\begin{lemma}\label{L:2}
For any $j\in A$, $T(A,\lambda,j)\cap S_A^\lambda\subseteq S'_{A\setmin\set{j}}$.
\end{lemma}

\begin{proof}
The proof is straightforward.
\end{proof}

We also put
$U(A,\lambda,i)=\Co\bigl(\set{p_i}\cup\setm{p(i,A,j)}{j\in A,j\ne i}\bigr)$.

\begin{lemma}\label{L:3}
For any $i\in A$,
$U(A,\lambda,i)\subseteq\bigcap\setm{T(A,\lambda,j)}{j\in A,j\ne i}$.
\end{lemma}

\begin{proof}
For any $j\in A$, $j\ne i$, the polytope $T(A,\lambda,j)$ contains the point
$p_i$ and the point $p(i,A,j)$. Moreover, it contains the whole face
$S_{A\setmin\set{j}}$ whence all the points $p(i,A,k)$, $k\ne i,j$. Therefore,
$U(A,\lambda,i)\subseteq T(A,\lambda,j)$, for all $j\in A$, $j\ne i$.
\end{proof}

\begin{lemma}\label{L:4}
For any $i,j\in A$ such that $i\ne j$,
$U(A,\lambda,i)\cap S'_{A\setmin\set{j}}=\set{p(i,A,j)}$.
\end{lemma}

\begin{proof}
$p(i,A,j)\in U(A,\lambda,i)\cap S'_{A\setmin\set{j}}$ by the definition
of $U(A,\lambda,i)$ and
$S'_{A\setmin\set{j}}$. To prove the reverse inclusion, we suppose that
$z\in U(A,\lambda,i)\cap S'_{A\setmin\set{j}}$. Then there are
$\mu_j\in[0,1]$, $j\in A$, such that
$\sum\setm{\mu_j}{j\in A}=1$ and
$z=\mu_i p_i+\sum\setm{\mu_j p(i,A,j)}{j\in A,j\ne i}$. Since
$S'_{A\setmin\set{j}}$ is
a face and $p_i\notin S'_{A\setmin\set{j}}$, we have $\mu_i=0$ and
\[
\setm{p(i,A,j)}{j\in A,j\ne i,\mu_j\ne 0}\subseteq S'_{A\setmin\set{j}}.
\]
Obviously, $p(i,A,k)\notin S'_{A\setmin\set{j}}$, for all $k\ne i,j$. Thus,
$\mu_k=0$, for all $k\ne i,j$, whence $\mu_j=1$ and $z=p(i,A,j)$.
\end{proof}

\begin{lemma}\label{L:5}
If $q_i\in U(A,\lambda,i)\setmin\setm{p(i,A,j)}{j\in A,j\ne i}$, for all
$i\in A$, then $S_A^\lambda\subseteq\interior{|A|}\Co\bigl(\setm{q_i}{i\in A}\bigr)$.
\end{lemma}

\begin{proof}
For any $i\in A$, we put
$B_i=\Co\bigl(\setm{q_j}{j\in A,j\ne i}\bigr)$. Then $B_i\subseteq T(A,\lambda,i)$, for all
$i\in A$, by Lemma~\ref{L:4}. Moreover, if
$B_i\cap S'_{A\setmin\set{i}}\ne\varnothing$,
then there extsts $j\in A\setmin\set{i}$ such that
$q_j\in S'_{A\setmin\set{i}}\cap U(A,\lambda,j)$ since
$S'_{A\setmin\set{i}}$ is a face of $T(A,\lambda,i)$. By Lemma~\ref{L:4}, this
implies that $q_j=p(j,A,i)$, a contradiction with the choice of $q_j$. Therefore,
$B_i\subseteq T(A,\lambda,i)\setmin S'_{A\setmin\set{i}}$.

By Lemma~\ref{L:2}, we get $S_A^\lambda\cap B_i=\varnothing$, for all $i\in A$.
Thus, for any $i\in A$, $S_A^\lambda$ is a subset of the open half-space $X_i$
defined by the hyperplane which contains $B_i$. Hence,
$S_A^\lambda\subseteq\bigcap\setm{X_i}{i\in A}=\interior{|A|}\Co\bigl(\setm{q_i}{i\in A}\bigr)$.
\end{proof}

\begin{lemma}\label{L:6}
There is $\varepsilon(\lambda)>0$ such that
$S_A^\lambda\subseteq\interior{|A|}\Co\bigl(S_{A\setmin\set{i}}^\varepsilon
\cup S_{A\setmin\set{j}}^\varepsilon\bigr)$,
for any $\varepsilon\in(0,\varepsilon(\lambda)]$ and any $i,j\in A$, $i\ne j$.
\end{lemma}

\begin{proof}
We pick $\varepsilon(\lambda)>0$ with respect to the property that
the extreme point $p_k^{\varepsilon(\lambda)}$ of the polytope
$S_{A\setmin\set{i}}^{\varepsilon(\lambda)}$
(of the polytope $S_{A\setmin\set{j}}^{\varepsilon(\lambda)}$, respectively)
belongs to $U(A,\lambda,k)$, for all $k\in A\setmin\set{i}$
(for all $k\in A\setmin\set{j}$, respectively).
The desired conclusion follows then from Lemma~\ref{L:5}.
\end{proof}

We construct the finite set $X$ which provides an embedding of the lattice
$\Subm{n+1}$ into the lattice $\CoX$. Let $v$ be
the center of $S_{\mathbf{n+1}}$. Let $\lambda_0>0$ be small enough.
Suppose that $k<n-1$ and we have already found
$\lambda_0$,\dots, $\lambda_k>0$ such that
$\lambda_j\in(0,\varepsilon(\lambda_{j-1})]$, for all $0<j\leqslant k$. By
Lemma~\ref{L:6}, there exists $\lambda_{k+1}\in(0,\varepsilon(\lambda_k)]$
such that, for any $A\subseteq\mathbf{n+1}$ with $|A|=n+1-k>2$ and any
$i,j\in A$, $i\ne j$, we have
$S_A^{\lambda_k}\subseteq\interior{|A|}\Co\bigl(S_{A\setmin\set{i}}^{\lambda_{k+1}}
\cup S_{A\setmin\set{j}}^{\lambda_{k+1}}\bigr)$.
We put $\lambda_n=0$. For any nonempty $A\subseteq\mathbf{n+1}$ and any $i\in A$,
we also put
\[
P_A=S_A^{\lambda_k},\qquad U(A,i)=U(A,\lambda_k,i),\qquad p(i,A)=p_i^{\lambda_k}
\]
where $k<n+1$ is such that $|A|+k=n+1$.

\begin{lemma}\label{L:7}
For any $A\subseteq B\subseteq\mathbf{n+1}$ and any $i\in A$, we have
$U(A,i)\subseteq U(B,i)$.
\end{lemma}

\begin{proof}
We argue by induction on $|B\setmin A|$. If $|B\setmin A|=0$ then
$U(B,i)=U(A,i)$, and we are done. Let $j\in B\setmin A$. By the induction
hypothesis, $U(A,i)\subseteq U(B\setmin\set{j},i)$. All the extreme points
of the polytope $U(B\setmin\set{j},i)$ are in the interior of the face of
$U(B,i)$ which is the convex hull of the set
$\set{p_i}\cup\setm{p(i,B,k)}{k\in B,k\ne i,j}$. Therefore,
$U(B\setmin\set{j},i)\subseteq U(B,i)$.
\end{proof}

We define the desired set $X$ by
\[
X=\set{v}\cup\bigcup\setm{\Ex P_A}{A\subset\mathbf{n+1}}.
\]

First we notice the important property of the lattice $\CoX$.

We remind that the \emph{join dependency relation} $\DD$ is defined
for join irreducible elements $a,b$ of a lattice $L$, $a\DD b$, if $a \not = b$,
and there is a $p \in L$ with
$a \leq b \vee p$ and $a \not \leq c \vee p$ for $c < p$. 
A $\DD$-sequence is a finite sequence $a_0,\dots,a_{n-1}$ $(n\geq 2)$
of join irreducible
elements of $L$ such that $a_i\DD a_{i+1}$ for all $i <n$, where the subscripts
are computed modulo $n$.
It is well-known
that a finite lattice $L$ is lower bounded iff it contains no $\DD$-cycles
(see, for example, Corollary 2.39 in \cite{FJN}).
 
\begin{lemma}\label{L:LB}
The finite lattice $\CoX$ is lower bounded.
\end{lemma}

\begin{proof}
If $a,b\in X\setmin\set{v}$, then there are $A,B\subseteq\mathbf{n+1}$ such that
$a\in\Ex P_A$ and $b\in\Ex P_B$. In this case,
$\set{a}\DD\set{b}$ implies that $|B|<|A|$. Moreover, $\set{v}\DD\set{a}$,
for any $a\in X\setmin\set{v}$, and $\set{a}\DD\set{v}$ holds for no
$a\in X$. Thus, the lattice $\CoX$ does not contain a $\DD$-cycle whence it is
lower bounded.
\end{proof}

Secondly, we observe that the composition of $\psi_X$ defined in section 2, and
$\varphi$ given by (2) is a a desired mapping of lattices.

\begin{proposition}\label{step2}
The map $\psi_X\varphi\colon\Subm{n+1}\to\CoX$ is a lattice embedding.
\end{proposition}

\begin{proof}
Since both $\psi_X$ and $\varphi$ preserve meets, the composition
$\psi_X\varphi$ also does.

If $A\in B_0\setmin B_1$, for some $B_0,B_1\in\Subm{n+1}$, then
$x\in\psi_X\varphi(B_0)\setmin\psi_X\varphi(B_1)$, where
$x\in\Ex P_{\mathbf{n+1}\setmin A}$
in the case $A\subset\mathbf{n+1}$ and $x=v$ in the case $A=\mathbf{n+1}$.
Therefore, the map $\psi_X\varphi$ is one-to-one.

To prove that $\psi_X\varphi$ preserves joins, it suffices to show that, for any
noncomparable sets $A_0,A_1\subseteq\mathbf{n+1}$,
\[
\psi(A_0\cap A_1)\cap X\subseteq\Co\bigl(\psi(A_0)\cup\psi(A_1)\bigr)\cap X,
\]
where $\psi$ is the map defined by~\eqref{psi}.
By the definition, we have
\[
\psi(A_0\cap A_1)\cap X=\Ex P_{A_0\cup A_1}=\setm{p(i,A_0\cup A_1)}{i\in
A_0\cup A_1},
\]
when $A_0\cup A_1\subset\mathbf{n+1}$, and
\[
\psi(A_0\cap A_1)\cap X=\set{v},
\]
when $A_0\cup A_1=\mathbf{n+1}$.
By Lemma~\ref{L:7}, for any $j_i\in A_i$, $i<2$, we have
$p(j_i,A_i)\in U(A_i\cup\set{j_{1-i}},j_i)\subseteq U(A_0\cup A_1,j_i)$.
Thus, by Lemma~\ref{L:5}, we get
\begin{align*}
\psi(A_0\cap A_1)\cap X&\subseteq
\Co\bigl(\setm{p(i,A_0)}{i\in A_0}\cup\setm{p(i,A_1)}{i\in A_1}\bigr)\cap X\\
&=\Co\bigl(\psi(A_0)\cup\psi(A_1)\bigr)\cap X.
\end{align*}
Moreover, for any $A_0,A_1\subseteq\mathbf{n+1}$ such that
$A_0\cup A_1=\mathbf{n+1}$, we have that
$v\in\Co\bigl(\psi(A_0)\cup\psi(A_1)\bigr)$. The proof of the lemma is
complete.
\end{proof}

Now we state the main result of this section.

\begin{theorem}
For any finite lower bounded lattice $L$, there is $n\in\omega$ and
a finite set $X\subseteq\R^n$ such that the lattice $\Co(\R^n,X)$ is lower
bounded and $L$ embeds into both $\Co(\R^n)$ and $\Co(\R^n,X)$.
\end{theorem}

\begin{proof}
According to~\cite{Ad,Rep},
for any finite lower bounded lattice $L$, there is
$n\in\omega$ such that $L$ is isomorphic to a sublattice of
$\Subm{n+1}$. The desired conclusion follows from
Propositoins~\ref{step1} and~\ref{step2}.
\end{proof}

{\bf Acknowledgments.} The author wants to thank F.~Wehrung for the question that
he sent in June of 2000 that inspired the construction of
Proposition~\ref{step1}, also for the
consequent discussion of the idea of proof of Proposition~\ref{step2}.
We are also grateful
to G.~Bergman whose wonderful paper~\cite{B} sparkled the fruitful communication
on the topic and motivated the writing up of these results.
Many valuable suggestions about the reorganization of proofs were sent to us by
M.~Semenova. They are implemented in the current version of the paper.

\end{document}